\documentclass{amsart}
\usepackage{amsmath,amscd,amsopn,amsthm,amsxtra,upref,amssymb,latexsym}
\theoremstyle{plain}
\newtheorem{theorem}{Theorem}

\newtheorem{lemma}{Lemma}

\theoremstyle{definition}

\theoremstyle{remark}

\numberwithin{equation}{section}
\newcommand{\s}{\sigma}
\newcommand{\R}{\mathbb R}
\newcommand{\Rn}{\mathbb R^n}
\newcommand{\RM}{\mathbb R^{n+1}_+}
\newcommand{\Rm}{\mathbb R^{n+1}}
\newcommand{\al}{\alpha}
\begin{document}
\title{Doubling properties of caloric functions}
\author{L. Escauriaza}
\address[L. Escauriaza]{Universidad del Pa{\'\i}s Vasco / Euskal Herriko
Unibertsitatea\\Dpto. Matematicas\\Apto. 644, 48080 Bilbao, Spain.} 
\email{mtpeszul@lg.ehu.es}
\thanks{The first and second authors are supported  by MEC grant MTM2004-03029 and by the
European Commission via the network Harmonic Analysis and Related Problems, project number
RTN2-2001-00315. The third author is supported by MIUR grant 2002013279.}
\author{F.J. Fern\'andez}
\address[F.J. Fern\'andez]{Ceremade\\Universite de Paris IX- Dauphine\\UMR CNRS 7534, France.} 
\email{fernandez@ceremade.dauphine.fr}
\author{S. Vessella}
\address[S. Vessella]{DiMaD\\Universit\`a di Firenze\\Via Lombroso 6/17, 50134 Firenze,
Italy.} 
\email{sergio.vessella@dmd.unifi.it}
\dedicatory{"Science is the great antidote to the poison of enthusiasm and
superstition" (Adam Smith, Wealth of Nations, 5, III, 3)}
\begin{abstract}
    We obtain quantitative estimates of unique continuation for solutions to
parabolic equations: doubling properties and two-sphere one-cylinder inequalities.
\end{abstract}
\maketitle

\section{Introduction}\label{S:1}
This work is concerned with  quantitative estimates of unique continuation for  parabolic equations. In the context
of elliptic equations the quantitative estimates of unique continuation, which have been derived from the nowadays
standard methods to prove qualitative results of unique
continuation ({\it Carleman} or {\it frequency function}  methods)  are the so called
\emph{doubling properties} and 
\emph{optimal three spheres inequalities}. In particular, it is by now well known that when
$$E=\sum_{i,j=1}^n\partial_i(g^{ij}(x)\partial_j\ )$$ is an elliptic operator in $\Rn$, there is a constant $N$
depending on $n$, $M$, the parameters of ellipticity and the Lipschitz norm of the matrix of coefficients of the
operator $E$ such that, solutions of the inequality, $|Eu|\le M\left(|u|+|\nabla u|\right)$ on $B_4$,
verify the
\emph{doubling property} \cite{gl} 
\begin{equation*}\int_{B_{2r}(z)}u^2\,dx\le (N\Theta)^N \int_{B_{r}(z)}u^2\,dx\ ,\end{equation*}
when $|z|\le 1$, $0<r\le 1/2$ and 
\begin{equation*}\Theta
={\int_{B_4}u^2\,dx}/{\int_{B_1}u^2\,dx}\end{equation*}
and the \emph{optimal three sphere inequality} (\cite{lo}, \cite{l})
\begin{equation*}
\int_{B_1}u^2\,dx\le \left(\int_{B_r}u^2\,dx\right)^{\frac 1{1+N\log{\left(1/r\right)}}}
\left(N\int_{B_4}u^2\,dX\right)^{\frac{N\log{\left(1/r\right)}}{1+N\log{\left(1/r\right)}}}\ ,
\end{equation*}
when $0<r\le 1$.

The three sphere inequality is called optimal for the following reason. If one seeks to find the largest
possible function $\theta :(0,1)\longrightarrow (0,1)$ such that the three sphere inequality
\begin{equation}\label{E1}
\int_{B_1}u^2\,dx\le \left(\int_{B_r}u^2\,dx\right)^{\theta(r)}
\left(N\int_{B_4}u^2\, dX\right)^{1-\theta(r)}
\end{equation}
holds for some positive constant $N$, for all harmonic functions $u$ in $B_4$ and $0<r<1$,
one finds from the identity
\begin{equation*}\int_{B_r}u^2(x)\,dx=\frac {r^{2k+n}}{2k+n}\int_{\partial B_1}u^2\,d\sigma\ 
,\end{equation*} which holds when $u$ is a homogeneous harmonic polynomial of degree $k\ge 1$, and \eqref{E1}
that
\begin{equation*}
1/4\le
\left( r/4\right)^\theta N^{\frac{1-\theta}{2k+n}}\quad\text{for all}\quad k\ge 1\ . 
\end{equation*}
Thus, $\theta(r)$ cannot be larger than $\log
4/\log{(4/r)}$. On the other hand, the recent applications to stability issues in Inverse Problems of
quantitative estimates of unique continuation (\cite{abrv}, \cite{amr}) show it is useful to know that \lq\lq
harmonic\rq\rq functions are as small as they can be at  intermediate scales when  they are known to be small at
smaller scales.

The first quantitative result of {\it strong} unique continuation for parabolic equations,  a {\it
two-sphere one-cylinder inequality} \cite[Theorem $9'$]{lo} (see the body of the paper for the relevant
definitions),  was derived in 1974 from the corresponding first qualitative results of {\it
strong} unique continuation for second order parabolic equations
 in the literature \cite{lo}. In \cite{lo}, E.M. Landis and O.A. Oleinik used a reduction argument to the already
established elliptic results of unique continuation, which relays on the representation formula of solutions to
parabolic equations with {\it time independent} coefficients in terms of the eigenfunctions of the
corresponding elliptic operator, and thus their results did not apply to  parabolic equations with {\it
time dependent} coefficients. After 1988 when more results of unique continuation for
elliptic equations were available in the literature, the qualitative results
in \cite{lo} were reproved with less regularity assumptions in 
\cite{l}. The quantitative results in \cite{lo} were also reprove with less regularity
assumptions in 
\cite{crv}. The authors in
\cite{l} and
\cite{crv} used the same reduction argument as in \cite{lo} and as a consequence, their results only apply
to parabolic operators with {\it time independent} coefficients.\par
In the same way as it has not been a trivial task to derive quantitative
estimates of unique continuation for elliptic
operators from the corresponding qualitative results (See \cite{gl}
and \cite{l}), the same happens in the parabolic case.

In \cite{p}, \cite{c}, \cite{e}, \cite{ev}, \cite{ef}, \cite{f} and \cite{av} the authors have obtained with
standard arguments of unique continuation, {\it frequency functions or Carleman methods}, qualitative results of
{\it strong} unique continuation for parabolic equations with {\it time dependent} coefficients. In
particular, between \cite{ef} and \cite{f} the following qualitative-quantitative result is
proved in relation with backward parabolic operators
\begin{equation*}
P=\sum_{i,j=1}^n\partial_i\left(g^{ij}(X)\partial_j\ \right)+\partial_t\ ,
\end{equation*}
which for some positive constants $\lambda$ and $M$ satisfy the conditions
\begin{gather}
\lambda
|\xi|^2\le\sum_{i,j=1}^ng^{ij}(X)\xi_i\xi_j\le\lambda^{-1}|\xi|^2\ ,\label{E:parabolicidad}\\
\sum_{i,j=1}^n|g^{ij}(X)-g^{ij}(Y)|\le
M\left(|x-y|^2+|t-s|\right)^{1/2}\ ,\label{E:regularidad}
\end{gather}
when $X=(x,t)$, $Y=(y,s)$ are in $\mathbb R^{n+1}$ and $\xi$ is in $\Rn$.

\begin{theorem}\label{T:qualitative} Let $P$ be a backward second order parabolic operator as
above. Then, there is a constant $N=N(\lambda , M, n)$ such that, if $u$ satisfies
$$|Pu|\le M\left(|u|+|\nabla u|\right)$$
in $Q_2$ and $u(x,0)$ has a zero of infinite order at $x=0$, the following holds in $Q_1$
$$|u(x,t)|\le Ne^{-1/(Nt)}\|u\|_{L^\infty(Q_2)}\quad .$$  
\end{theorem}

Here and in the sequel $B_r(z)=\{\ x\in\Rn\ :\ |x-z|<r\ \}$, $B_r=B_r(0)$,
$Q_r(z, \tau)=B_r(z)\times [\tau,\tau+r^2]$ and
$Q_r=Q_r(0,0)$.

Theorem \ref{T:qualitative} gives a
qualitative-quantitative result of strong unique continuation  but is missing the proper and natural
quantitative estimates  of strong unique continuation for
second order parabolic equations: {\it doubling properties} within characteristic hyperplanes and {\it two-sphere
one-cylinder inequalities}. Here, we prove these quantitative estimates with new
arguments, which are based on the {\it frequency functions or Carleman inequalities}
appearing in
\cite{p},
 \cite{ef} and \cite{f}  and extend their applicability to solutions of parabolic
equations with {\it time dependent} coefficients. In particular, the following theorem is proved:

\begin{theorem}\label{T:quantitative} Let $P$ be a backward second order parabolic operator
verifying \eqref{E:parabolicidad} and \eqref{E:regularidad} and $u$ satisfy
$|Pu|\le M\left(|u|+|\nabla u|\right)$ in $Q_4$.
Then, there is a constant $N= N(\lambda, M, n)$ such that the
following properties  hold when $0<r\le 1/2$
\begin{enumerate}
\item $\int_{B_{2r}}u^2(x,0)\,dx\le \left(N\Theta\right)^N\int_{B_{r}}u^2(x,0)\,dx$.
\item $\int_{B_1}u^2(x,0)\,dx\le \left(\int_{B_r}u^2(x,0)\,dx\right)^{\frac 1{1+N\log{\left(1/r\right)}}}
\left(N\int_{Q_4}u^2\,dX\right)^{\frac{N\log{\left(1/r\right)}}{1+N\log{\left(1/r\right)}}}$.
\item If $u|_{t=0}$ is not identically zero, then $|u|_{|_{t=0}}$ is a Muckenhoupt
weight in $B_{1/2}$.
\end{enumerate}
\par
Here, $\Theta=\int_{Q_4}u^2\,dX/\int_{B_1}u^2(x,0)\,dx$ and $dX=dxdt$.
\end{theorem}

Recall that a  nonzero and locally integrable function $w$ in $B_1$  is called a Muckenhoupt
weight when there are numbers 
$p\in (1,+\infty)$ and $C>0$ such that
\begin{equation*}
\left(\text{\rlap |{$\int_{B_r(z)}$}}w^p\,dx\right)^{\frac 1p}\le C
\text{\rlap |{$\int_{B_r(z)}$}}w\ dx\quad \text{and}\quad \int_{B_{2r(z)}}w\,dx\le C\int_{B_r(z)}w\,dx
\end{equation*}
for all balls $B_r(z)$ such that $B_{4r}(z)\subset B_1$. Property 2 generalizes the analogue result for elliptic
operators in \cite{gl}.\par
In regard to {\it space-time} like {\it doubling properties} we have the following result.
\begin{theorem}\label{T:space-time doubling} Assume that $P$ and $u$ satisfy the conditions in Theorem
\ref{T:quantitative}. Then, there is $N= N(\lambda, M, n)$ such that the
following  holds when $0<r\le 1/\sqrt{N\log{\left(N\Theta\right)}}$
\begin{enumerate}
\item $\int_{Q_{r}}u^2\,dX\le
e^{N\log{\left(N\Theta\right)}\log{\left(N\log{\left(N\Theta\right)}\right)}}r^2\int_{B_{r}}u^2(x,0)\,dx$.
\item $\int_{Q_{2r}}u^2\,dX\le
e^{N\log{\left(N\Theta\right)}\log{\left(N\log{\left(N\Theta\right)}\right)}}\int_{Q_{r}}u^2\,dX$.
\end{enumerate}
\end{theorem}
We do not know if Theorem \ref{T:space-time doubling} is optimal. In fact, if $u$ is a backward
caloric polynomial of degree $k\ge 1$, (i.e. $u(\lambda x,\lambda^2t)=\lambda^ku(x,t)$ when
$(x,t)$ is in
$\Rm$, $\lambda>0$) the following happens when $r>0$
\begin{equation*}
\frac{\int_{Q_{r}}u^2\,dX}{r^2\int_{B_{r}}u^2(x,0)\,dx}=\frac{\int_{Q_{1}}u^2\,dX}{\int_{B_{1}}u^2(x,0)\,dx}
\quad .
\end{equation*}\par
Theorems \ref{T:quantitative} and \ref{T:space-time doubling}
also hold when $u$ verifies $|Pu|\le M\left(|u|+|\nabla u|\right)$ in $Q_4\cap
D\times [0,+\infty)$, $D=\{\ (x',x_n) : x_n>\varphi(x')\ \}$, where
$\varphi:\R^{n-1}\longrightarrow\R$ is a $C^{1,1}$ function verifying $\varphi(0)=0$ and
when either the Dirichlet or Neumann (conormal derivative) data of $u$ vanishes
identically on $\partial D\times [0,+\infty)\cap Q_4$.  In fact, when the matrix of
coefficients of $P$ is time independent and $u$ has zero Dirichlet condition on the
lateral boundary the results hold when $\varphi=\varphi_1+\varphi_2$, where $\varphi_1$
is a convex function and $\varphi_2$ is $C^{1,\al}$ for some $\al$ in (0,1). Of course,
one must replace $B_r$ by $B_r\cap D$ and so on in the corresponding statement of the
boundary version of the Theorems. We do no prove these results here but they
follow combining the arguments in this work with others in
\cite{ef}.
\par  In section
\ref{S:2} we prove Theorems
\ref{T:quantitative} and
\ref{T:space-time doubling} when $P=\Delta +\partial_t$ with the frequency function
method. With some extra work it is possible to prove Theorems \ref{T:quantitative} and
\ref{T:space-time doubling} with the same method and when the main coefficients of $P$ are
time independent. In section
\ref{S:3}, the more general case is proved using a more detailed version of the Carleman
inequalities appearing in
\cite{ef} and \cite{f}.\par 
\section{Constant coefficients}\label{S:2} 
 The proof of Theorems \ref{T:quantitative} and \ref{T:space-time doubling}  when
$P=\Delta+\partial_t$ is a consequence of the following four Lemmas.
\begin{lemma}\label{L:decay}
Assume that $u$ satisfies $\left |\Delta u+\partial_tu\right |\le M \left(|u|+|\nabla u|\right)$
in $Q_4$. Then, there is $N=N(n,M)$ such that the inequalities
\begin{equation*}
N\log(N\Theta)\ge 1\quad\text{and}\quad N\int_{B_2}u^2(x,t)\,dx\ge\int_{B_1}u^2(x,0)\,dx
\end{equation*}
hold when $ t\le
1/\left(N\log(N\Theta)\right)$ and 
$\Theta=\int_{Q_4}u^2\,dX/\int_{B_1}u^2(x,0)\,dx$.
\end{lemma}
\begin{proof}
 Let $f=u{\varphi}$, where
$\varphi\in C_0^\infty(B_2)$,
$0\le\varphi
\le 1$,
$\varphi=1$ in
$B_{3/2}$. The standard estimate \cite{ls}
\begin{equation}\label{E:standard estimate}\|u\|_{L^\infty(Q_{3})}+\|\nabla
u\|_{L^\infty(Q_{3})}\le N\|u\|_{L^2(Q_{4})}
\end{equation}
 for solutions  of the
inequality
\begin{equation}\label{E:heat}
\left |\Delta u+\partial_tu\right |\le M \left(|u|+|\nabla u|\right)\ \text{in}\ Q_4
\end{equation}
gives
\begin{equation}\label{E:heat2}
|\Delta f+\partial_tf|\le N\left[ |f|+|\nabla f|+\|u\|_{L^2
(Q_4)}\chi_{B_2\setminus B_{3/2}}\right]
\end{equation}
and shows that the first claim holds.
Setting $H(t)=\int f^2(x,t)G(x-y,t)\,dx$,
where
$G(x,t)=t^{-n/2}e^{-|x|^2/4t}$ and
$y\in B_1$ we have
\begin{equation}\label{E:heat3}\dot H(t)=2\int f(\Delta f+\partial_tf)G(x-y,t)\,dx+2\int |\nabla
f|^2G(x-y,t)\,dx\end{equation}
and from \eqref{E:heat2}, \eqref{E:heat3} and the Cauchy-Schwarz's inequality
\begin{equation*}\dot H(t)\ge -N\left(H(t)+e^{-1/Nt}\|u\|_{L^2(Q_4)}^2\right)\quad .\end{equation*}
Integration of this inequality in $(0,t)$ for $t\in (0,16)$ gives
\begin{equation*}
N\int f^2(x,t)G(x-y,t)\,dx\ge u^2(y,0)-Ne^{-1/Nt}\|u\|_{L^2(Q_4)}^2\quad .
\end{equation*}
Integrating this inequality over $B_1$ and recalling that $\int
G(x-y,t)\,dy=1$ we get
\begin{equation*}
N\int_{B_2}u^2(x,t)\,dx\ge \int_{B_1}u^2(x,0)\,dx-Ne^{-1/Nt}\|u\|_{L^2(Q_4)}^2\quad ,
\end{equation*}
and the Lemma follows from this inequality and \eqref{E:standard estimate}.
\end{proof}
\begin{lemma}\label{L:freq}
Given $a>0$ and $f\in W^{2,\infty}(\RM)$ set $G_a=\left(t+a\right)^{-n/2}e^{-|x|^2/4\left(t+a\right)}$,
\begin{equation*}
H_a(t)=\int_{\Rn}f^2G_a\,dx\ ,\ D_a(t)=\int_{\Rn}|\nabla f|^2G_a\,dx
\ \text{and}\  N_a(t)=\frac{2(t+a)D_a(t)}{H_a(t)}\ .
\end{equation*}
Then, 
\begin{equation*}
\dot N_a(t)\ge -\frac{(t+a)}{H_a(t)}\int (\Delta
f+\partial_tf)^2G_a\, dx\quad .
\end{equation*}
\end{lemma}
\begin{proof}
The identities $\partial_tG_a-\Delta G_a=0$, $\nabla G_a=-\frac x{2\left(t+a\right)}G_a$,
$\Delta=\text{div}\left(\nabla\ \right)$ and integration by parts imply the following
identities
\begin{equation}\label{E:deriv}
\dot H_a(t)=2\int f(\Delta f+\partial_tf)G_adx+2D_a(t)\ , 
\end{equation}
\begin{align*}\dot
H_a(t)&=2\int f\left(\partial_tf+\tfrac x{2(t+a)}\cdot
\nabla f-\tfrac 12\left(\Delta f+\partial_tf\right)\right)G_a\,dx+\int f\left(\Delta f+\partial_t
f\right)G_a\,dx\ ,\\ D_a(t)&=\int f\left(\partial_tf+\tfrac {x}{2(t+a)}\cdot \nabla
f-\tfrac 12 \left(\Delta f+\partial_tf\right)\right)G_a\,dx-\tfrac 12\int f\left(\Delta f+\partial_t
f\right)G_a\,dx
\end{align*}
and
\begin{align}\label{E:mult}
\dot H_a(t)D_a(t)&=2\left(\int f\left(\partial_tf+\tfrac x{2(t+a)}\cdot \nabla f-\tfrac 12\left(\Delta
f+\partial_tf\right)\right)G_a\,dx\right)^2\\&-\tfrac 12\left(\int f\left(\Delta f+\partial_t
f\right)G_a\,dx\right)^2\ .\notag
\end{align}
The Rellich-N\v{e}cas identity with vector field $\nabla G_a$ 
\begin{align*}
\text{div}(\nabla G_a|\nabla f|^2)&-2\text{div}((\nabla f\cdot\nabla G_a)\nabla
f)\\&=\Delta G_a|\nabla f|^2-2D^2G_a\nabla f\cdot\nabla f-2\nabla f\cdot\nabla G_a\Delta f
\end{align*}
and integration by parts give
\begin{align}\label{E:Rellich}
\int\Delta G_a&|\nabla f|^2\,dx=2\int D^2G_a\nabla f\cdot\nabla f\,dx+2\int\nabla
f\cdot\nabla G_a\Delta f\,dx\\&=2\int \left(\tfrac x{2(t+a)}\cdot\nabla
f\right)^2G_a\,dx-2\int\tfrac x{2(t+a)}\cdot\nabla f\Delta fG_a\,dx-D_a(t)/(t+a)\ .\notag
\end{align}
Again, the fact that $G_a$ is a caloric function, integration by parts, \eqref{E:Rellich} and the completion of
the square of $\partial_tf+\tfrac {x}{2(t+a)}\cdot
\nabla f-\tfrac 12\left(\Delta f+\partial_tf\right)$ yields the formula

\begin{align}\label{E:derv2}
\dot D_a(t) &=2\int\left(\partial_tf+\tfrac
{x}{2(t+a)}\cdot
\nabla f-\tfrac 12\left(\Delta f +\partial_tf\right)\right)^2G_a\,dx\\&-\tfrac 12\int\left(\Delta
f+\partial_tf\right)^2G_a\,dx-D_a(t)/(t+a)\ .\notag
\end{align}
Then, from \eqref{E:mult},\eqref{E:derv2} and the quotient rule
\begin{align}\label{derv5}
\dot N_a(t)&=\frac {4(t+a)}{H_a(t)^2}\Biggl\{\int\left(\partial_tf+\tfrac
{x}{2(t+a)}\cdot\nabla f-\tfrac 12 \left(\Delta
f+\partial_tf\right)\right)^2G_a\,dxH_a(t)\\ &-\left(\int f\left(\partial_tf+\tfrac
{x}{2(t+a)}\cdot \nabla f-\tfrac 12\left(\Delta
f+\partial_tf\right)\right)G_a\,dx\right)^2\notag\\& +\tfrac{1}{4}\left(\int f\left(\Delta
f+\partial_tf\right)G_a\,dx\right)^2-\tfrac{1}4\int\left(\Delta
f+\partial_tf\right)^2G_a\,dxH_a(t)\Biggr\}\notag
\end{align}
and Lemma \ref{L:freq} follows from \eqref{derv5}, the Cauchy-Schwarz inequality and the
positiveness of the third term on the right hand side of \eqref{derv5}.
\end{proof}
\begin{lemma}\label{L:ineq} The inequality
\begin{equation*}
\int \tfrac {|x|^2}{8a}h^2e^{-|x|^2/4a}\,dx\le 2a\int |\nabla
h|^2e^{-|x|^2/4a}\,dx+\tfrac n2\int h^2e^{-|x|^2/4a}\,dx
\end{equation*}
holds for all $h\in C_0^\infty(\Rn)$ and $a>0$.
\end{lemma}
 \begin{proof}
 The inequality follows setting  $v=he^{-|x|^2/8a}$ and from the identity
\begin{align*}
2a\int |\nabla
h|^2e^{-|x|^2/4a}\,dx+&\tfrac n2\int h^2e^{-|x|^2/4a}\,dx-\int \tfrac
{|x|^2}{8a}h^2e^{-|x|^2/4a}\,dx\\&=2a\int |\nabla v|^2\,dx\ .
\end{align*}
\end{proof}
\begin{lemma}\label{L:implidoub} Assume that $N$ and $\Theta$ verify
$N\log(N\Theta)\ge 1$,
$h\in C_0^\infty(\Rn)$ and that the inequality 
\begin{equation*}
\begin{split}
2a\int |\nabla h|^2e^{-|x|^2/4a}\,dx&+\tfrac n2\int
h^2e^{-|x|^2/4a}\,dx\le N\log{(N\Theta)}\int
h^2e^{-|x|^2/4a}\,dx
\end{split}
\end{equation*}
holds when $a\le\frac 1{12N\log{(N\Theta)}}$. Then,
\begin{equation*}
\int_{B_{2r}}h^2\,dx\le (N\Theta)^N\int_{B_r}h^2\,dx\ \text{when}\ 0<r\le 1/2\ .
\end{equation*}
\end{lemma}
\begin{proof}
The inequality satisfied by $h$ and Lemma \ref{L:ineq} show that  
\begin{equation*}
\int \tfrac{|x|^2}{8a}h^2e^{-|x|^2/4a}\,dx\le N\log{(N\Theta)}\int
h^2e^{-|x|^2/4a}\,dx
\end{equation*}
when $a\le1/\left(12N\log{(N\Theta)}\right)$.
For given $0<r\le 1/2$ and $0< a\le
\tfrac{r^2}{16N\log{\left(N\Theta\right)}}\
$, the last inequality  implies 
\begin{gather*}
\int\tfrac{|x|^2}{8a}h^2e^{-|x|^2/4a}\,dx\le N\log{(N\Theta)}\left[\int_{B_r}
h^2\,dx+\tfrac{8a}{r^2}\int_{\Rn\setminus B_r}
\tfrac{|x|^2}{8a}h^2e^{-|x|^2/4a}\,dx\right]\\\le N\log{(N\Theta)}\int_{B_r} h^2\,dx+\tfrac
12\int
\tfrac{|x|^2}{8a}h^2e^{-|x|^2/4a}\,dx\ .
\end{gather*}
Thus,
\begin{equation}\label{E:casiadob}
\int \tfrac{|x|^2}{16a}h^2e^{-|x|^2/4a}\,dx\le N\log{(N\Theta)}\int_{B_r}
h^2\,dx\quad \text{when}\quad  0< a\le
\tfrac{r^2}{16N\log{\left(N\Theta\right)}}\ .
\end{equation}
Now, $e^{-|x|^2/4a}|x|^2/(16a)\ge
(N\Theta)^{-N}N\log{(N\Theta)}$ when $r\le |x|\le 2r$ and $a=
\tfrac{r^2}{16N\log{\left(N\Theta\right)}}\ $. This and \eqref{E:casiadob}
imply
\begin{equation*}
\int_{B_{2r}}h^2\,dx\le (N\Theta)^N\int_{B_r}h^2\,dx\ \text{when}\ 0<r\le 1/2\ .
\end{equation*}
\end{proof}
\begin{proof}[Proof of Theorems \ref{T:quantitative} and \ref{T:space-time
doubling}. Constant coefficients] Let $u$ satisfy
\eqref{E:heat} and set
$f=u\psi$ in Lemma
\ref{L:freq}, where
$\psi\in C_0^\infty(B_4)$,
$0\le\psi\le 1$, $\psi=1$ in $B_3$ and $\psi=0$ outside $B_{\frac 72}$. Then,
\begin{equation}\label{E:heat4}
|\Delta f+\partial_tf|\le N\left[ |f|+|\nabla f|+\|u\|_{L^2
(Q_4)}\chi_{B_4\setminus B_{3}}\right]
\end{equation}
and from  Lemmas \ref{L:decay},\  \ref{L:freq} and \eqref{E:heat4}, we have
\begin{equation}\label{E:nece}
\begin{split}
H_a(t)&\ge N^{-1}
(t+a)^{-n/2}e^{-1/(t+a)}\int_{B_1}u^2(x,0)\,dx\ ,\\
\dot N_a(t)&\ge -N\left[1+N_a(t)+e^{-5/4(t+a)}\Theta\right]\ ,
\end{split}
\end{equation}
 when $0\le t+a\le 1/N\log(N\Theta)$ and where $\Theta$ was defined  in Lemma \ref{L:decay}.
The choice in Lemma \ref{L:decay} of the constant $N$ gives,
$e^{-5/4(t+a)}\Theta\le 1$ when
$0\le t+a\le 1/N\log(N\Theta)$. This fact and the second inequality in \eqref{E:nece} imply that
\begin{equation}\label{E:incres}
e^{Nt}N_a(t)+Ne^{Nt}\ \text{is non-decreasing
when}\ t+a\le 1/N\log(N\Theta)\quad .
\end{equation}
Multiplying \eqref{E:deriv} by $\frac{t+a}{H_a(t)}\ $, we find from \eqref{E:heat4}, Lemma
\ref{L:decay} and the first inequality in \eqref{E:nece} that there is some constant $N=N(n,M)$ such that
\begin{equation}\label{E:relacion}
N_a(t)\le N\left[1+(t+a)\partial_t\log {H_a(t)}\right]\ ,\ \text{when}\ \ 0\le t+a\le
1/N\log(N\Theta)
\end{equation}
and
\begin{equation}\label{E:desigualdad importante}
\partial_t\log {H_a(t)}\le N\left[1+N_a(t)\right]/(t+a)\ ,\ \text{when}\ \ 0\le t+a\le
1/N\log(N\Theta)\quad .
\end{equation}

Setting $\beta =\frac 1{N\log{(N\Theta )}}\ $, we get from \eqref{E:incres}, \eqref{E:relacion} and the
first inequality in \eqref{E:nece} that 
\begin{equation}\label{E:acotacion frecuencia}
\begin{split}
N_a(t)&\lesssim N_a(\beta/4)+1\lesssim 1+\int_{\beta /4}^{\beta
/2}\frac{N_a(t)}{(t+a)}\,dt\lesssim 1+\int_{\beta /4}^{\beta
/2}\partial_t\log{H_a(t)}\,dt\\&=1+\log{\left(\tfrac{H_a(\beta/2)}{H_a(\beta/4)}\right)}
\le N\log{(N\Theta)}\ ,\ \text{when}\quad t+a\le\beta/12
\end{split}
\end{equation}
and with constants depending only on $M$ and $n$. In particular, the
inequality
\begin{equation}\label{E:medsob}
\begin{split}
2a\int |\nabla f(x,0)|^2e^{-|x|^2/4a}\,dx&+\tfrac n2\int
f^2(x,0)e^{-|x|^2/4a}\,dx\\&\le N\log{(N\Theta)}\int
f^2(x,0)e^{-|x|^2/4a}\,dx
\end{split}
\end{equation}
holds for all  $0<a\le 1/\left(12N\log{(N\Theta)}\right)$

From Lemma \ref{L:implidoub} with $h=f(\ \cdot,0)$ and  recalling that $f(x,0)=u(x,0)$ in
$B_3$ we obtain the doubling property claimed in Theorem \ref{T:quantitative} when $z=0$
\begin{equation}\label{E:doubling}
\int_{B_{2r}}u^2(x,0)\,dx\le (N\Theta)^N\int_{B_r}u^2(x,0)\,dx\ \text{when}\ 0<r\le 1/2\quad  .
\end{equation}

For the same values of $r$
choose
$k\ge 2$ such that
$2^{-k}< r\le 2^{-k+1}$ and iterate \eqref{E:doubling} when $r=2^{-j}$, $j=0, \dots ,k-1$. It gives 
\begin{equation*}
\int_{B_1}u^2(x,0)\,dx\le (N\Theta)^{N\log{\left(1/r\right)}}\int_{B_r}u^2(x,0)\,dx\quad ,
\end{equation*}
but writing the value of $\Theta$ in the above inequality, which was defined in Lemma \ref{L:decay},  one
gets 
\begin{equation*}
\int_{B_1}u^2(x,0)\,dx\le \left(\int_{B_r}u^2(x,0)\,dx\right)^{\frac 1{1+N\log{\left(1/r\right)}}}
\left(N\int_{Q_4}u^2\,dX\right)^{\frac{N\log{\left(1/r\right)}}{1+N\log{\left(1/r\right)}}}\quad,
\end{equation*}
which proves the two-sphere one-cylinder inequality.

Rewrite the integral on the left hand
side of the following inequalities as the sum of the corresponding integrals over
$B_r$ and 
$\Rn\setminus B_r$ and assume that $a\le r^2/\left(16N\log{\left(N\Theta\right)}\right)$, then 
\begin{gather*}
\begin{split}
\int f^2(x,0)e^{-|x|^2/4a}\,dx\le&\int_{B_r}f^2(x,0)\,dx+\tfrac {16a}{r^2}\int_{\Rn\setminus B_r}
\tfrac{|x|^2}{16a}f^2(x,0)e^{-|x|^2/4a}\,dx\\\le& \int_{B_r}f^2(x,0)\,dx+\tfrac
{1}{N\log{\left(N\Theta\right)}}\int
\tfrac{|x|^2}{16a}f^2(x,0)e^{-|x|^2/4a}\,dx
\end{split}
\end{gather*}
and recalling that \eqref{E:casiadob} holds when $h=f(\ \cdot,0)$ and
$a\le r^2/\left(16N\log{\left(N\Theta\right)}\right)$, it follows that
\begin{equation}\label{E:importante}
\int f^2(x,0)e^{-|x|^2/4a}\,dx\le 2\int_{B_r}u^2(x,0)\,dx\quad \text{when}\quad  0< a\le
\tfrac{r^2}{16N\log{\left(N\Theta\right)}}\ .
\end{equation}

The fact
that $f(x,0)=u(x,0)$ when $|x|\le 2$ and choosing $a=
\tfrac{r^2}{16N\log{\left(N\Theta\right)}}$ in the last inequality and \eqref{E:medsob}   
implies 
\begin{equation}\label{E:cacciopoli}
\int_{B_r}|\nabla u(x,0)|^2\,dx\le (N\Theta)^Nr^{-2}\int_{B_r}u^2(x,0)\,dx\quad .
\end{equation}
Then, \eqref{E:cacciopoli} and the Sobolev inequality
\begin{equation*}
\left(\text{\rlap |{$\int_{B_r}$}}|\varphi-\text{\rlap
|{$\int_{B_r}$}}\varphi|^{\frac{2n}{n-2}}\,dx\right)^{\frac{n-2}{2n}}\le Nr\left(\text{\rlap
|{$\int_{B_r}$}}|\nabla \varphi |^2d\,x\right)^{\frac 12}\ ,\ \varphi\in C^\infty(B_r)
\end{equation*}
give
\begin{equation*}
\left(\text{\rlap |{$\int_{B_r}$}}|u(x,0)|^{\frac{2n}{n-2}}\,dx\right)^{\frac{n-2}{n}}\le
(N\Theta)^N\text{\rlap |{$\int_{B_r}$}}u^2(x,0)\,dx\ ,
\end{equation*}
which shows that $|u|^2|_{t=0}$ satisfies a reverse H\"older inequality centered at $0$. When $n=2$,
replace the previous Sobolev inequality by 
\begin{equation*}
\left(\text{\rlap |{$\int_{B_r}$}}|\varphi-\text{\rlap
|{$\int_{B_r}$}}\varphi|^6\,dx\right)^{\frac{1}{6}}\le Nr\left(\text{\rlap
|{$\int_{B_r}$}}|\nabla \varphi |^2d\,x\right)^{\frac 12}\ ,\ \varphi\in C^\infty(B_r)\quad . 
\end{equation*}

These  and other standard arguments show that $|u|^2|_{t=0}$ is a
Muckenhoupt weight  in $B_{1/2}$, which satisfies a reverse h\"older inequality from
$L^{\frac n{n-2}}$-averages to $L^1$-averages and with constant 
\begin{equation*}
\left(\frac{N\int_{Q_4}u^2\,dX}{\int_{B_{1/2}}u^2(x,0)\,dx}\right)^N\quad .
\end{equation*}
To
prove Theorem
\ref{T:space-time doubling},
\eqref{E:acotacion frecuencia} and  
\eqref{E:desigualdad importante} imply
\begin{equation*}
\partial_t\log {H_a(t)}\le N\log(N\Theta)/(t+a)\ ,\ \text{when}\ \ 0< t+a\le
1/N\log(N\Theta)
\end{equation*}
and integration of this inequality over $[0,t]$
gives
\begin{equation*}
\int f^2(x,t)e^{-|x|^2/4(a+t)}\,dx\le \left(1+\tfrac
ta\right)^{N\log{\left(N\Theta\right)}}\int
f^2(x,0)e^{-|x|^2/4a}\,dx\ ,
\end{equation*}
when $0<t+a\le \frac 1{12N\log{\left(N\Theta\right)}}\ $. From \eqref{E:importante} and the last inequality 
\begin{equation}\label{E: a integrar}
\int_{B_{2r}} u^2(x,t)\,dx\le \left(N\Theta\right)^N\left(1+\tfrac
ta\right)^{N\log{\left(N\Theta\right)}}\int_{B_r}u^2(x,0)\,dx\ ,
\end{equation}
when $0<t<4r^2\le \frac 1{32N\log{\left(N\Theta\right)}}\ $ and $a=
\tfrac{r^2}{16N\log{\left(N\Theta\right)}}$. The integration of \eqref{E: a integrar} over $[0,4r^2]$ gives
\begin{equation*}
\int_{Q_{2r}} u^2\,dX\le
\left(N\Theta\right)^Na\int_{B_r}u^2(x,0)\,dx\int_0^{4r^2/a}\left(1+s\right)^{N\log{\left(N\Theta\right)}}\,ds\ 
.
\end{equation*}
Thus,
\begin{equation}\label{E:relacion espacio-plano}
\int_{Q_{2r}}u^2\,dX\le 
e^{N\log{\left(N\Theta\right)}\log{\left(N\log{\left(N\Theta\right)}\right)}}r^2\int_{B_{r}}u^2(x,0)\,dx\ 
.
\end{equation}

The second part in Theorem \ref{T:space-time doubling} is a consequence of \eqref{E:relacion
espacio-plano}, the doubling property in Theorem  \ref{T:quantitative}  and the standard local
$L^{\infty}$-bounds of solutions of
\eqref{E:heat} in terms of their
$L^2$-averages \cite{ls}.
\end{proof}
\section{Variable coefficients}\label{S:3} 
In this section $P$ is a backward parabolic operator verifying
\eqref{E:parabolicidad} and
\eqref{E:regularidad} and $u$  a solution of the inequality
\begin{equation}\label{E:heatvariable}
|Pu|\le M \left(|u|+|\nabla u|\right)\ \text{in}\ Q_4\quad .
\end{equation}\par

As in section \ref{S:2} we need a few Lemmas to prove Theorems
\ref{T:quantitative} and \ref{T:space-time doubling}. The first one, Lemma
\ref{L:decayvariable}, is the variable coefficients version of Lemma
\ref{L:decay}. The second, Lemma \ref{L:Carlemaninequalty}, is a more detailed version
of the Carleman inequalities established in \cite{ef} and \cite{f}. These two imply Lemma
\ref{L:medsob2}, a result analogous to
\eqref{E:medsob}, which as we saw in section \ref{S:2},  implies Theorems
\ref{T:quantitative} and \ref{T:space-time doubling}.
\begin{lemma}\label{L:decayvariable}
There is $N=N(\lambda, M, n)$ such
that the inequalities
\begin{equation*}
N\log(N\Omega_\rho)\ge 1\quad\text{and}\quad
N\int_{B_{2\rho}}u^2(x,t)\,dx\ge\int_{B_\rho}u^2(x,0)\,dx
\end{equation*}
hold when $ 0<t\le
\rho^2/\left(N\log(N\Omega_\rho)\right)$ and $0<\rho\le 1$. Here,
\begin{equation*}
\Omega_\rho=\frac{\int_{Q_4}u^2\,dX}{\rho^2\int_{B_\rho}u^2(x,0)\,dx}\ .
\end{equation*}
\end{lemma}
\begin{proof}
	When $\rho=1$ set as in Lemma \ref{L:decay} $f=u{\varphi}$, where $\varphi\in
C_0^\infty(B_2)$,
$0\le\varphi
\le 1$ and
$\varphi=1$ in
$B_{3/2}$. Define $$H(t)=\int f^2(x,t)G(x,t;y,0)\,dx\ ,$$
where now
$G(x,t;y,s)$ is the fundamental solution in $\Rm$ of the parabolic operator $P^*$, the adjoint
of $P$,
$y\in B_1$ and $P^*G(x,t;y,s)=-\delta_{(y,s)}$ is the Dirac operator at $(y,s)$.\par The 
Gaussian bounds for the fundamental solution \cite{a}
\begin{equation*}
N^{-1}(t-s)^{-\frac n2}e^{-\frac{N|x-y|^2}{t-s}}\le G(x,t;y,s)\le N
(t-s)^{-\frac n2}e^{-\frac{|x-y|^2}{N(t-s)}}
\end{equation*}
and the fact that \eqref{E:standard estimate} is  satisfied in this setting \cite{ls}, show
that the arguments in the proof of Lemma
\ref{L:decay} can be repeated again, which proves Lemma
\ref{L:decayvariable} when $\rho=1$. The case $0<\rho\le 1$ follows  rescaling to the case
$\rho=1$ and from the observation
\begin{equation*}
\int_{Q_{4\rho}}u^2\,dX\le\int_{Q_{4}}u^2\,dX\ .
\end{equation*}
\end{proof}

In Lemma \ref{L:Carlemaninequalty} we use the following notation: for a given function $\s$ defined on
some interval and
$a>0$,
$\s_a(t)=\sigma(t+a)$, $G(x,t)=t^{-n/2}e^{-|x|^2/4t}$ is the Gauss kernel and  $G_a(x,t)=G(x,t+a)$.
\begin{lemma}\label{L:Carlemaninequalty}
Assume $g^{ij}(0,0)=\delta_{ij}$. Then, there is a constant $N=N(\lambda, M, n)$ verifying
the following property:  for each number
$\al\ge 2$  there is a  nondecreasing function
$\sigma :(0,+\infty)\longrightarrow\R_+$ satisfying $t/N\le \sigma (t)\le t$
when
$0<\al t\le 4$ and such that the inequality

\begin{gather*}
\al^2\int\s_a^{-\al}v^2G_a\,dX+\al\int \s_a^{1-\al}\ |\nabla
v|^2G_a\,dX\\\le N\int\s_a^{1-\al}|Pv|^2G_a\,dX+ N^{\al}\al^{\al}\sup_{t\ge 0}\int
v^2+|\nabla v|^2\,dx\\+\s(a)^{-\al}\left[-(a/N)\int|\nabla
v(x,0)|^2G(x,a)\,dx+N\al\int v ^2(x,0)G(x,a)\,dx\right]
\end{gather*}
holds for all $0 < a\le\frac 1{\al}$ and $v\in C_0^\infty(\Rn\times[0,\frac 3\al))$.
\end{lemma}
The proof of  Lemma \ref{L:Carlemaninequalty} is sketched at the end of this section.
We now prove  Lemma
\ref{L:medsob2}, where as in \eqref{E:medsob}, $f=u\psi$ with $\psi\in C_0^\infty(B_4)$,
$0\le\psi\le 1$, $\psi=1$ in $B_3$ and $\psi=0$ outside $B_{\frac 72}$.
\begin{lemma}\label{L:medsob2} There are constants $N= N(\lambda, M, n)\ge 1$  and
$\rho=\rho(\lambda, M, n)$ in (0,1) such that
\begin{equation*}
\begin{split}
2a\int |\nabla f(x,0)|^2e^{-|x|^2/4a}\,dx&+\tfrac n2\int
f^2(x,0)e^{-|x|^2/4a}\,dx\\&\le N\log{(N\Theta_\rho)}\int
f^2(x,0)e^{-|x|^2/4a}\,dx\ ,
\end{split}
\end{equation*}
 when $0< a\le 1/\left(12N\log{(N\Theta_\rho)}\right)$ and
$N\log{\left(N\Theta_\rho\right)}\ge 1$. Here,
\begin{equation*}
\Theta_\rho=\int_{Q_4}u^2\,dX/\int_{B_\rho}u^2(x,0)\,dx\ .
\end{equation*}
\end{lemma}
\begin{proof} For
fixed
$\al\ge 2$ and
$a\le \frac 1\al\
$ take as $v$ in the Carleman inequality  the function
$v=f(x,t)\varphi(t)=u(x,t)\psi(x)\varphi(t)$, where $\varphi=1$ when $t\le\frac 1\al$ and $\varphi=0$
when $t\ge \frac 2\al\ $. From \eqref{E:heatvariable} and \eqref{E:standard estimate} we have  
\begin{equation}\label{E:error}
|Pv|\le N\left[|v|+|\nabla v|+\|u\|_{L^2
(Q_4)}\chi_{B_4\times [0,\frac 2\al]\setminus B_{3}\times [0,\frac 1\al]}\right]\ . 
\end{equation}
Then, \eqref{E:error}, the observation that $\s^{1-\al}_aG_a\le N^{\al+\frac n2}\al^{\al+\frac n2-1}$ in the
region
$B_4\times [0,\frac 2\al]\setminus B_{3}\times [0,\frac 1\al]$, standard arguments used
with Carleman inequalities and the fact that $t/N\le \sigma (t)\le t$ on $(0,\frac 4\al)$,
imply that the inequality
\begin{equation}\label{E:applicarleman}
\begin{split}
&\al^2\int_0^{\frac 1\al}\int_{B_2}\left(t+a\right)^{-\al}u^2e^{-|x|^2/4(t+a)}\,dX\le
N^{\al}\al^{\al}\|u\|_{L^2 (Q_4)}^2\\&+\s (a)^{-\al}\left [-(a/N)\int |\nabla
f(x,0)|^2e^{-|x|^2/4a}\,dx+N\al\int f ^2(x,0)e^{-|x|^2/4a}\,dx\right]
\end{split}
\end{equation}
holds when
$\al\ge \frac n2 +2$ and $0<a\le \frac 1\al\ $.\par
For $0<\rho\le 1$, which will be chosen later, let
$N$ be a fixed constant larger than the ones appearing in Lemma
\ref{L:decayvariable} and \eqref{E:applicarleman} and assume that $\al\ge
N\log{(N\Omega_\rho)}$ and
$0<a\le
\frac{\rho^2}{2\al}$. Then,
\begin{equation}\label{E:applicarleman22}
\begin{split}
\al^2\int_0^{\frac 1\al}\int_{B_2}\left(t+a\right)^{-\al}u^2e^{-|x|^2/4(t+a)}\,dX&\ge
\al^2\int_0^{\frac
{\rho^2}{\al}}\int_{B_{2\rho}}\left(t+a\right)^{-\al}e^{-\frac{\rho^2}{(t+a)}}u^2\,dX\\
\ge \tfrac{\al^2}N\int_a^{a+\frac
{\rho^2}{\al}}s^{-\al}e^{-\frac{\rho^2}{s}}\,ds\int_{B_{\rho}}u^2(x,0)\,dx&\ge
\tfrac{\al^2}{N}\int_{\frac {\rho^2}{2\al}}^{\frac
{\rho^2}{\al}}s^{-\al}e^{-\frac{\rho^2}{s}}\,ds\int_{B_\rho}u^2(x,0)\,dx\\&\ge\frac{\al^{\al+1}}{2N}
\left(\frac
1{\rho e}\right)^{2\al}\frac{\|u\|_{L^2
(Q_4)}^2}{\Omega_\rho}\ .  
\end{split}
\end{equation}
The inequalities \eqref{E:applicarleman} and \eqref {E:applicarleman22} show that  to make
sure that the left hand side of
\eqref{E:applicarleman} is greater than  four times the  first term on right hand side of
\eqref{E:applicarleman} when $\al\ge N\log{(N\Omega_\rho)}$ and $0<a\le
\frac{\rho^2}{2\al}$, suffices to know that
\begin{equation}\label{E:suffices}
\left(\frac 1{\rho^2}\right)^{\al}\ge 8N\left(Ne^2\right)^\al\Omega_\rho\ .   
\end{equation}
Choose then $\rho$ as the solution of the equation $\frac 1{\rho^2}=8Ne^2$. Then,
\eqref{E:suffices}
 holds when $8^{\al-1}\ge N\Omega_\rho$. Thus, there are fixed constants,
$\rho=\rho(\lambda, M, n)$ in $(0,1)$ and $N=N(\lambda, M, n)\ge 1$, such that
\begin{equation*}
\begin{split}
&\tfrac 12\int_0^{\frac 1\al}\int_{B_2}\left(t+a\right)^{-\al}u^2e^{-|x|^2/4(t+a)}\,dX+
N^{\al}\al^{\al}\|u\|_{L^2 (Q_4)}^2\\&\le \s (a)^{-\al}\left [-(a/N)\int |\nabla
f(x,0)|^2e^{-|x|^2/4a}\,dx+N\al\int f ^2(x,0)e^{-|x|^2/4a}\,dx\right]\ ,
\end{split}
\end{equation*}
when $\al\ge N\log{(N\Omega_\rho)}$ and $0<a\le \frac{\rho^2}{12\al}\ $. In particular,
the inequality
\begin{equation}\label{E:applicarleman3}
\begin{split}
2a\int |\nabla f(x,0)|^2e^{-|x|^2/4a}\,dx&+\tfrac n2\int
f^2(x,0)e^{-|x|^2/4a}\,dx\\&\le N\log{(N\Omega_\rho)}\int
f^2(x,0)e^{-|x|^2/4a}\,dx
\end{split}
\end{equation}
is true when $0< a\le \rho^2/\left(12N\log{(N\Omega_\rho)}\right)$ and
$N\log{\left(N\Omega_\rho\right)}\ge 1$. Recall that
$\Omega_\rho=\Theta_\rho/\rho^2$ and rename the fraction $N/\rho^2$ as $N$ in the above
claim. Then, Lemma \ref{L:medsob2} follows.
\end{proof}
\begin{proof}[Proof of Theorems \ref{T:quantitative} and \ref{T:space-time
doubling}. Variable coefficients]
Lemmas \ref{L:implidoub} and \ref{L:medsob2} imply that 
\begin{equation}\label{E:doubling2}
\int_{B_{2r}}u^2(x,0)\,dx\le (N\Theta_\rho)^N\int_{B_r}u^2(x,0)\,dx\ \text{when}\ 0<r\le
1/2
\end{equation}
and \eqref{E:doubling2} and the  argument used in section \ref{S:2} to prove the
interpolation inequality or two-sphere one-cylinder inequality (2) in Theorem
\ref{T:quantitative}, show that the interpolation inequality
\begin{equation}\label{E:desigualdad de interpolacion}
\int_{B_\rho}u^2(x,0)\,dx\le \left(\int_{B_r}u^2(x,0)\,dx\right)^{\frac
1{1+N\log{\left(\rho/r\right)}}}
\left(N\int_{Q_4}u^2\,dX\right)^{\frac{N\log{\left(\rho/r\right)}}{1+N\log{
\left(\rho/r\right)}}}
\end{equation}   
holds when $r\le\frac\rho 2\ $. On the other hand, translations of $u$ in the space
variable show that the same interpolation inequality is true when we replace in
\eqref{E:desigualdad de interpolacion} the balls
$B_\rho$ and
$B_r$ by $B_\rho(x)$ and
$B_r(x)$ respectively, and where $x$ is any point in $B_1$.  Then, standard  covering
arguments combined with these interpolation inequalities at large scales
(\eqref{E:desigualdad de interpolacion} and its translated analogues  when
$r=\rho/2$), show that we can find  $\beta$ in $(0,1)$ such that
\begin{equation*}
\int_{B_1}u^2(x,0)\,dx\le \left(\int_{B_\rho}u^2(x,0)\,dx\right)^{\beta}
\left(N\int_{Q_4}u^2\,dX\right)^{1-\beta}\ ,
\end{equation*}
inequality which can be rewritten as $\left(N\Theta_\rho\right)^\beta\le N\Theta$. This
fact,
\eqref{E:doubling2} and the arguments in section \ref{S:2} finish the prove of
Theorem \ref{T:quantitative}. \par
Following the notation in Lemma \ref{L:freq}, the arguments in section \ref{S:2} give that
the claims 3 and 4 in Theorem
\ref{T:space-time
doubling} follow if we know that
\begin{equation}\label{E: a integrar despues}
\partial_t\log {H_a(t)}\le N\log(N\Theta)/(t+a)\ \ \text{when}\ \ 0< t+a\le
1/N\log(N\Theta)\ .
\end{equation}
On the other hand, translations of $u$ in the time variable, Lemma \ref{L:decayvariable}
and  Lemma \ref{L:medsob2} (In particular, applying to u(x,t+s) the 
arguments leading to the proof of \eqref{E:applicarleman3}, which satisfies
\eqref{E:heatvariable} for some backward parabolic operator $P$ when $s\le
\rho^2/N\log(N\Omega_\rho)$, then recalling that in this
case $f(x,0)=u(x,s)\psi(x)$ and then, rewriting the variables $a$ and $s$ as $a+t$ and $t$
respectively), show that with a possibly a larger 
$N$ and a smaller $\rho$ 
\begin{equation*}
(a+t)\int |\nabla
f(x,t)|^2e^{-|x|^2/4(a+t)}\,dx\le N\log{\left(N\Omega_\rho\right)}\int f
^2(x,t)e^{-|x|^2/4(a+t)}\,dx
\end{equation*}
when $0< t+a\le
\rho^2/N\log(N\Omega_\rho)$. The last  inequality,  Lemma \ref{L:ineq} and the inequality
$\left(N\Theta_\rho\right)^\beta\le N\Theta$ imply that
\begin{equation}\label{E:acotacion frecuencia2}
\begin{split}
N_a(t)&\le N\log{\left(N\Theta\right)}\ ,\\
\int\tfrac {|x|^2}{a+t}f^2(x,t)G_a\,dx&\le  N\log(N\Theta)\int f^2(x,t)G_a\,dx\ ,
\end{split}
\end{equation}
when $0< t+a\le
1/N\log(N\Theta)$.\par
Only using the Lipschitz regularity in the space-variable of the matrix of coefficients
of
$P$ to bound the first order term arising in the formal calculation of $P^*G_a$, gives
that
$|P^*G_a|\le N(|x|^2+t+a)/(t+a)^2$. This fact,
\eqref{E:acotacion frecuencia2} and a calculation similar to the one leading to
\eqref{E:deriv}, which was done to compute $\dot H_a(t)$ in terms of $D_a(t)$, but which on
this occasion is done to bound
$\dot H_a(t)$  in terms of
$D_a(t)$  but  replacing the backward heat operator by
$P$, imply that \eqref{E: a integrar despues} holds and prove
 Theorem \ref{T:space-time
doubling}. 
\end{proof}
\begin{proof}[Proof of Lemma \ref{L:Carlemaninequalty}]
Here we do not give a complete proof of Lemma \ref{L:Carlemaninequalty} since it is
basically already done in
\cite{ef} and
\cite{f}. We only outline the main details.\par 
As in \cite{ef}, if $P$ is a backward parabolic
operator satisfying \eqref{E:parabolicidad},
\eqref{E:regularidad} and $g^{ij}(0,0)=\delta_{ij}$, we consider
the time independent backward parabolic operator
$Q$
\begin{equation*}
Q=\sum_{i,j=1}^n\partial_i(g^{ij}(x)\partial_j\ )+\partial_t\quad \text{where}\quad
g^{ij}(x)=g^{ij}(x,0)\ , \ x\in\Rn\ .
\end{equation*}

In \cite[Theorem 4]{ef} it is shown with an integration by parts argument that there are
$N= N(\lambda, M, n)$ and $\delta_0=\delta_0(\lambda, M,n)$ such that the inequality
\begin{equation}\label{E: Carleman primitiva}
\begin{split}
\delta^{-2}\int\s^{2-\al}|D^2v|^2G\,dX+\al&\int\s^{-\al}\tfrac{\theta(\gamma
t)}tv^2G\,dX+\int \s^{1-\al}\tfrac{\theta(\gamma t)}t|\nabla v|^2G\,dX\\\le
N\int\s^{1-\al}|Qv|^2G\,dX + & N^{\al}\gamma^{\al+N}\int
\left(v^2+t|\nabla v|^2\right)\,dX
\end{split}
\end{equation}
is satisfied when $v\in C_0^\infty(\Rn\times (0,\frac 3{\gamma}))$,
$\gamma=\al/\delta^2$, $\al\ge 2$,
$0<\delta\le \delta_0$ and where $\s$ is the solution of the ordinary differential equation
\begin{equation*}
\tfrac
d{dt}\left[\log{\left(\tfrac\s{t\dot\s}\right)}\right]=\tfrac{\theta(\gamma t)}t\ ,\ 
\s(0)=0\  ,\ \dot\s(0)=1\ ,\text{where}\ 
\theta(t)=t^{\frac 12}\left(\log{\tfrac 5t}\right)^{\frac 32}\ .
\end{equation*}
In fact, $\s(t)=\beta(\gamma t)/\gamma$ with
\begin{equation*}
\beta(t)=t\ \text{exp}\left[-\int_0^{
t}\left(1-\text{exp}\left(-\int_0^s\tfrac{\theta (u)}u\ du\right)\right)\tfrac
{ds}s\right]\ .
\end{equation*}
 The
main point here is that for some $N$, $t/N\le \beta (t)\le t$, when $0<t\le 4$ and
\begin{equation}\label{E: acotaci—n sigma}
t/N\le \s (t)\le t\ ,\ \text{when}\ 0<t\le 4/\gamma\ .
\end{equation} 

If for fixed
$0<a\le \tfrac 1{\gamma}$ one repeats the calculations which led to \eqref{E: Carleman
primitiva} in  \cite{ef}, but now working  with  $v\in C_0^\infty(\Rn\times [0,\frac
3{\gamma}))$ and carrying out the integration over
$\Rn\times [a,+\infty)$, one must take into account the boundary terms which occur when
integrating by parts with respect to the time-variable and add up on the
on the right hand side of \eqref{E: Carleman primitiva} certain new terms. In fact, an analysis 
of the proof of Theorem 4 in \cite{ef} shows that there is $N=N(\lambda, M, n)$, such that the
sum of those boundary terms is bounded from above by
\begin{gather}\label{E: terminos frontera en tiempo}
\s(a)^{-\al}\left[-(a/N)\int |\nabla
v(x,a)|^2G(x,a)\,dx+N\al \int v^2(x,a)G(x,a)\,dx
\right.\\\left. +N\int \tfrac{|x|^3}{a}\,v^2(x,a)G(x,a)\,dx\right]\ .\notag
\end{gather}
\eqref{E: acotaci—n sigma} implies, 
$\tfrac{|x|^3}{a}\,G(x,a)\s (a)^{-\alpha}\le N^{\alpha}\gamma^{\alpha+N}$, when $|x|\ge \delta$, $0<a<1$. This
fact, Lemma \ref{L:ineq} and the writing of the third integral in
\eqref{E: terminos frontera en tiempo} as the sum of the same integral over $B_\delta$ and $\Rn\setminus
B_\delta$ show that \eqref{E: terminos frontera en tiempo} can be
bounded by 
\begin{gather}\label{E: terminos frontera en tiempoI}
\s(a)^{-\al}\left[-(a/N)\int |\nabla
v(x,a)|^2G(x,a)\,dx+N\al \int v^2(x,a)G(x,a)\,dx\right]
\\+N^{\alpha}\gamma^{\alpha+N}\int v^2(x,a)\,dx\ ,\notag
\end{gather}
when $\delta$ is sufficiently small. Thus, there
is
$N= N(\lambda, M, n)$ and $\delta_0=\delta(M,\lambda , n)$ such that if $v\in
C_0^\infty(\Rn\times [0,\frac 3{\gamma}))$, $\al\ge 2$, $0<a\le \tfrac 1{\gamma}$
and
$\delta\le \delta_0$, then  
\begin{gather*}
\delta^{-2}\int\s^{2-\al}|D^2v|^2G\,dX+\al\int\s^{-\al}\tfrac{\theta(\gamma
t)}tv^2G\,dX+\int \s^{1-\al}\tfrac{\theta(\gamma t)}t|\nabla v|^2G\,dX\\\le
N\int\s^{1-\al}|Qv|^2G\,dX +  N^{\al}\gamma^{\al +N}\sup_{t\ge a}\int
v^2+t|\nabla v|^2\,dx\\
+\s (a)^{-\al}\left [-(a/N)\int |\nabla
v(x,a)|^2G(x,a)\,dx+N\al\int v^2(x,a)G(x,a)\,dx\right]\ ,
\end{gather*}
and where the integration in the integrals with respect to Lebesgue
measure $dX$ is done over $\Rn\times [a,+\infty)$.

If we plug in the last inequality the function
$v(x,t-a)$, consider the change of variables $t=s+a$ and observe that $\tfrac{\theta(\gamma t)}t\ge \gamma$ when
$0<\gamma t<4$, we get upon renaming the new variable
$s$ as
$t$ that the  inequality
\begin{gather*}
\delta^{-2}\int\s_a^{2-\al}|D^2v|^2G_a\,dX+\al^2\int\s_a^{-\al}v^2G_a\,dX+\al\int
\s_a^{1-\al}\ |\nabla
v|^2G_a\,dX\\\le N\int\s_a^{1-\al}|Q(v)|^2G_a\,dX+ N^{\al}\gamma^{\al +N}\sup_{t\ge 0}\int
v^2+(t+a)|\nabla v|^2\,dx\\
+\s (a)^{-\al}\left [-(a/N)\int |\nabla
v(x,0)|^2G(x,a)\,dx+N\al\int v^2(x,0)G(x,a)\,dx\right]
\quad 
\end{gather*}
is satisfied when $\al\ge 2$, $v\in
C_0^\infty(\Rn\times [0,\frac 3{\gamma}))$, $0<a\le \tfrac 1{\gamma}$,
$\delta\le \delta_0$, $\delta_0$ is sufficiently small, and where the integration  with respect to Lebesgue
measure $dX$ is done over $\Rn\times [0,+\infty)$.\par
On the other hand,  
$|g^{ij}(x,t)-g^{ij}(x)|\le M\sqrt t$. Then, choosing $\delta$ sufficiently small  it is
possible to replace $Q$ by
$P$ on the right hand side of the above inequality, which finishes the proof of Lemma
\ref{L:Carlemaninequalty}.
\end{proof}
\newpage

\end{document}